\documentclass[a4paper]{amsart}
\usepackage[margin=3cm]{geometry}
\usepackage{amssymb}
\usepackage{amsthm}
\usepackage{bm}
\usepackage{mathtools}
\usepackage[export]{adjustbox}

\newcommand{\diff}[2]{\frac{\mathrm{d} #1}{\mathrm{d} #2}}
\newcommand{\ddiff}[2]{\frac{\mathrm{d}^2 #1}{\mathrm{d} #2^2}}
\newcommand{\pfrac}[2]{\frac{\partial #1}{\partial #2}}
\newcommand{\ppfrac}[2]{\frac{\partial^2 #1}{\partial #2^2}}
\newcommand{\eref}[1]{(\ref{#1})}

\newtheorem{theorem}{Theorem}[section]
\renewenvironment{proof}[1][Proof]{\begin{trivlist}
\item[\hskip \labelsep {\bfseries #1}]}{\end{trivlist}}
\renewcommand{\qed}{\nobreak \ifvmode \relax \else
      \ifdim\lastskip<1.5em \hskip-\lastskip
      \hfill \fi \nobreak
      \vrule height0.5em width0.5em depth0em\fi}

\theoremstyle{remark}
\newtheorem{remark}{Remark}[section]

\title[Travelling wave solutions to the Keller--Segel model]{A geometric construction of travelling wave solutions to the Keller--Segel model}

\author[K.~Harley \and P.~van Heijster \and G.~J.~Pettet]{K.~Harley$^{\ast}$ \and P.~van Heijster$^{\ast}$ \and G.~J.~Pettet$^{\ast}$}

\thanks{$^{\ast}$Mathematical Sciences School, Queensland University of Technology, Brisbane, QLD 4000 Australia}

\begin{document}
\maketitle

\begin{abstract}
We study a version of the Keller--Segel model for bacterial chemotaxis, for which exact travelling wave solutions are explicitly known in the zero attractant diffusion limit.  Using geometric singular perturbation theory, we construct travelling wave solutions in the small diffusion case that converge to these exact solutions in the singular limit.  
\end{abstract}

\section{Introduction}
\label{sec:intro}

The Keller--Segel model \cite{Keller_Segel_71a, Keller_Segel_71b} is a very popular model for modelling cell migration in response to a chemical gradient, see for example \cite{Hillen_Painter_09, Tindall_Maini_Porter_Armitage_08} and references therein.  Because it has exact travelling wave solutions in the limit $D_u \to 0$ \cite{Feltham_Chaplain_00}, we are interested in the following particular version of the Keller--Segel model: 
\begin{equation}\label{eq:KS}\begin{aligned}
\pfrac{u}{t} &= D_u\ppfrac{u}{x} - Kw, \\
\pfrac{w}{t} &= D_w\ppfrac{w}{x} - \pfrac{}{x}\left(\frac{\chi w}{u}\pfrac{u}{x}\right),
\end{aligned}\end{equation}
with $u > 0$, $w \geq 0$, $x \in \mathbb{R}$, $t > 0$, $K, \chi > 0$, $D_{u,w} \geq 0$.  Here $u(x,t)$ is the concentration of the chemical or chemoattractant and $w(x,t)$ is the density of the migrating species.  In particular, we are interested in finding travelling wave solutions to \eref{eq:KS} in the case where both the diffusivities are small but of the same order: $0 \leq D_{u,w} \ll 1$.  

With $D_{u,w}$ small, \eref{eq:KS} is a singularly perturbed system; due to the advection (chemotactic) term we are unable to scale out the small parameters.  This makes \eref{eq:KS} amenable for analysis via geometric singular perturbation theory ({\sc gspt}) \cite{Jones_95,Kaper_99}, and we show that it supports travelling wave solutions.  In the limit $D_u \to 0$, these solutions agree with the exact solutions given in \cite{Feltham_Chaplain_00}.  

The background states of \eref{eq:KS} are $(u,w) = (u^{\ast},0)$, with $u^{\ast} \geq 0$ for physically relevant solutions.  We are interested in travelling wave solutions and so introduce a comoving frame $z = x - ct$ and \eref{eq:KS} becomes
\begin{equation}\label{eq:tw-KS}\begin{aligned}
-c \diff{u}{z} &= D_u \ddiff{u}{z} - Kw, \\
-c \diff{w}{z} &= D_w \ddiff{w}{z} - \diff{}{z} \left(\frac{\chi w}{u}\diff{u}{z} \right). 
\end{aligned}\end{equation}
Travelling wave solutions satisfy
\begin{equation}\label{eq:bc}
\lim_{z\to -\infty}{u(z)} = u_l, \quad \lim_{z\to\infty}{u(z)} = u_r > u_l, \quad \lim_{z\to \pm\infty}{w(z)} = 0. 
\end{equation}
Assuming $u_r > u_l$ implies $c > 0$; that is, we look for right-moving travelling waves.  

\subsection{An exact solution for $D_u = 0$}
\label{subsec:exact-sol}

As alluded to above, for $D_u = 0$ and $D_w < \chi$, \eref{eq:tw-KS} has exact solutions given by
\begin{equation}\label{eq:exact-sols}\begin{aligned}
u(z) &= \left[\sigma_2 + \sigma_1 \exp\left(-\frac{cz}{D_w}\right) \right]^{\frac{D_w}{D_w - \chi}}, \\
w(z) &= A\exp\left(-\frac{cz}{D_w}\right) \left[ \sigma_2 + \sigma_1 \exp\left(-\frac{cz}{D_w} \right) \right]^{\frac{\chi}{D_w - \chi}}, 
\end{aligned}\end{equation}
where
\[ \sigma_1 = \frac{AK(\chi - D_w)}{c^2}, \quad \sigma_2 = \frac{B(D_w - \chi)}{D_w}, \]
and $A$ and $B$ are constants of integration \cite{Feltham_Chaplain_00}.  By taking the limit of $u(z)$ in \eref{eq:exact-sols} as $z \to \mp\infty$ for fixed $0 < D_w < \chi$, we determine that $u_l = 0$ and 
\[ B = \frac{D_w}{D_w - \chi}u_r^{\frac{D_w - \chi}{D_w}}. \]
Consequently, we redefine
\begin{equation}\label{eq:newsigma2}
\sigma_2 = u_r^{\frac{D_w - \chi}{D_w}}.  
\end{equation}

\begin{remark}\label{rem:origin}
With $u_l = 0$, a travelling wave solution connects $(0,0)$ to $(u_r,0)$.  Although \eref{eq:KS} is not defined at $u = 0$, the solutions are still well behaved as $u \to 0$ since
\[ \lim_{z\to-\infty}{\frac{w}{u}} = \frac{c^2}{K(\chi - D_w)}.  \]
See also Section~\ref{subsec:persistence}.
\end{remark}

\subsection{Taking the limit as $D_w \to 0$}
\label{subsec:limit}

\begin{figure}[t]
\centering
\includegraphics{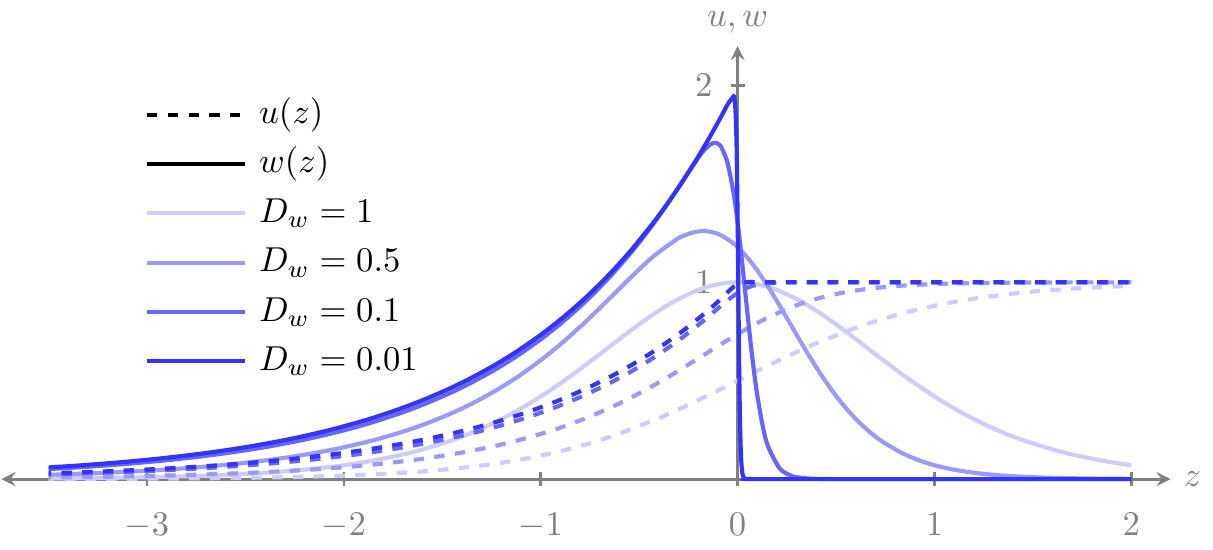}
\caption{Plots of $u(z)$ and $w(z)$ defined in \eref{eq:exact-sols}, with decreasing $D_w$ and parameters taken from \protect\cite{Feltham_Chaplain_00}: $\chi = 2$, $K = 1$, $c = 2$, $A = 4$, $u_r = 1$.}
\label{fig:shock-sol}
\end{figure}

Since we are interested in the case where both diffusivities are small, consider the limit of \eref{eq:exact-sols} as $D_w \to 0$.  Evaluating the limit gives 
\begin{equation}\label{eq:limit} 
\lim_{D_w \to 0}{u(z)} = \begin{cases} u_r{\rm e}^{cz/\chi}, & z \leq 0, \\ u_r, & z > 0, \end{cases} \quad {\rm and} \quad
\lim_{D_w \to 0}{w(z)} = \begin{cases} \dfrac{c^2 u_r}{K\chi}{\rm e}^{cz/\chi}, & z \leq 0, \\ 0, & z > 0, \end{cases} 
\end{equation}
which has a discontinuity or shock in $w$ at $z = 0$.  Figure~\ref{fig:shock-sol} shows solution curves of \eref{eq:exact-sols} for decreasing $D_w$, holding the other parameters constant.  

We now state our main result:
\begin{theorem}\label{th:existence}
Let $D_u = \mu\varepsilon$ and $D_w = \varepsilon$, with $0 < \varepsilon \ll 1$ a sufficiently small parameter and $\mu$ a positive, $\mathcal{O}(1)$ (with respect to $\varepsilon$) constant.  Then, travelling wave solutions to \eref{eq:KS} connecting $(0,0)$ to $(u_r,0)$ with $u_r > 0$, exist.  
\end{theorem}

\section{Geometric singular perturbation methods}

We use \textsc{gspt} to prove Theorem~\ref{th:existence}.  {\sc gspt} can be applied to problems exhibiting a clear separation of spatial scales; for example, cell migration where diffusion is operating on a much slower spatial scale than advection or reaction.  The power of this method lies in the ability to separate the spatial scales into independent, generically lower dimensional problems, which are more amenable to analysis.  

\begin{proof}
Following \cite{Wechselberger_Pettet_10}, we introduce a third variable $v = u_x$ such that
\begin{equation}\label{eq:bal-law}
\begin{pmatrix} u \\ v \\ w \end{pmatrix}_t + \begin{pmatrix} 0 \\ Kw \\ \chi vw/u \end{pmatrix}_x = \begin{pmatrix} -Kw \\ 0 \\ 0 \end{pmatrix} + \varepsilon \begin{pmatrix} \mu u \\ \mu v \\ w \end{pmatrix}_{xx}.  
\end{equation}
In the travelling wave coordinate, \eref{eq:bal-law} becomes
\[ \begin{aligned} 
\left(\mu\varepsilon u_z + cu\right)_z &= Kw, \\
\left(\mu\varepsilon v_z + cv - Kw\right)_z &= 0, \\
\left(\varepsilon w_z + cw - \frac{\chi vw}{u}\right)_z &= 0. \\
\end{aligned} \] 
The above system can be written as a system of first order differential equations by introducing the slow variables
\[ \begin{aligned}
\tilde{u} &\coloneqq \mu\varepsilon u_z + cu, \\
\tilde{v} &\coloneqq \mu\varepsilon v_z + cv - Kw, \\
\tilde{w} &\coloneqq \varepsilon w_z + cw - \frac{\chi vw}{u}, \\
\end{aligned} \]
to give
\[ \begin{aligned}
\mu\varepsilon u_z &= \tilde{u} - cu, \\
\mu\varepsilon v_z &= \tilde{v} - cv + Kw, \\
\varepsilon w_z &= \tilde{w} - cw + \frac{\chi vw}{u}, 
\end{aligned} \qquad \begin{aligned}
\tilde{u}_z &= Kw, \\
\tilde{v}_z &= 0, \\
\tilde{w}_z &= 0. 
\end{aligned} \]

The last two equations imply $\tilde{v}$ and $\tilde{w}$ are constants, which can be shown to be identically zero.  Thus, effectively we have a four-dimensional {\em slow system} in the slow travelling wave coordinate $z$:
\begin{equation}\label{eq:slow}
\begin{aligned}
\mu\varepsilon u_z &= \tilde{u} - cu, \\
\mu\varepsilon v_z &= -cv + Kw, \\
\varepsilon w_z &= -cw + \frac{\chi vw}{u}, \\
\tilde{u}_z &= Kw.   
\end{aligned}
\end{equation}
Equivalently, written in terms of the fast travelling wave coordinate $y = z/\varepsilon$ ($\varepsilon \neq 0$) we have the {\em fast system}:  
\begin{equation}\label{eq:fast}
\begin{aligned}
\mu u_y &= \tilde{u} - cu, \\
\mu v_y &= -cv + Kw, \\
w_y &= -cw + \frac{\chi vw}{u}, \\
\tilde{u}_y &= \varepsilon Kw. 
\end{aligned}
\end{equation}

In the singular limit the slow system reduces to 
\begin{equation}\label{eq:reduced} \begin{aligned}
0 &= \tilde{u} - cu, \\
0 &= -cv + Kw, \\
0 &= -cw + \frac{\chi vw}{u}, \\
\tilde{u}_z &= Kw,
\end{aligned} \end{equation}
which we call the {\em reduced problem}, and the fast system in the singular limit becomes
\begin{equation}\label{eq:layer} \begin{aligned}
\mu u_y &= \tilde{u} - cu, \\
\mu v_y &= -cv + Kw, \\
w_y &= -cw + \frac{\chi vw}{u}, \\
\tilde{u}_y &= 0, 
\end{aligned} \end{equation}
which we refer to as the {\em layer problem}.  Note that in the singular limit the two systems are no longer equivalent.  

\subsection{Layer problem}
\label{subsec:layer}

The steady states of the layer problem \eref{eq:layer} define a one-dimensional critical manifold $S$:
\begin{equation}\label{eq:S}
S = \left\{ (u,v,w,\tilde{u}) \, \middle| \, u = \frac{\tilde{u}}{c}, v = \frac{Kw}{c}, 0 = w\left(\frac{\chi v}{u} - c\right) \right\}, 
\end{equation}
where $\tilde{u}$ acts as a parameter.  This critical manifold has two distinct branches, 
\[ S_a \coloneqq \left\{ (u,v,w,\tilde{u}) \, \middle| \, u = \frac{\tilde{u}}{c}, v = 0, w = 0 \right\} \]
and 
\[ S_r \coloneqq \left\{ (u,v,w,\tilde{u}) \, \middle| \, u = \frac{\tilde{u}}{c}, v = \frac{\tilde{u}}{\chi}, w = \frac{c \tilde{u}}{\chi K} \right\}, \]
which intersect at $(u,v,w,\tilde{u}) = (0,0,0,0)$.  By examining the eigenvalues of the Jacobian of the linearised system, we can determine that $S_r$ is repelling, while $S_a$ is attracting, hence the subscript choice.  Thus, for each $\tilde{u}$, the layer flow connects a point on $S_r$ to the corresponding point on $S_a$, along what is referred to as a {\em fast fibre}.  
\begin{figure}[ht]
\centering
\includegraphics[width=0.49\textwidth,valign=t]{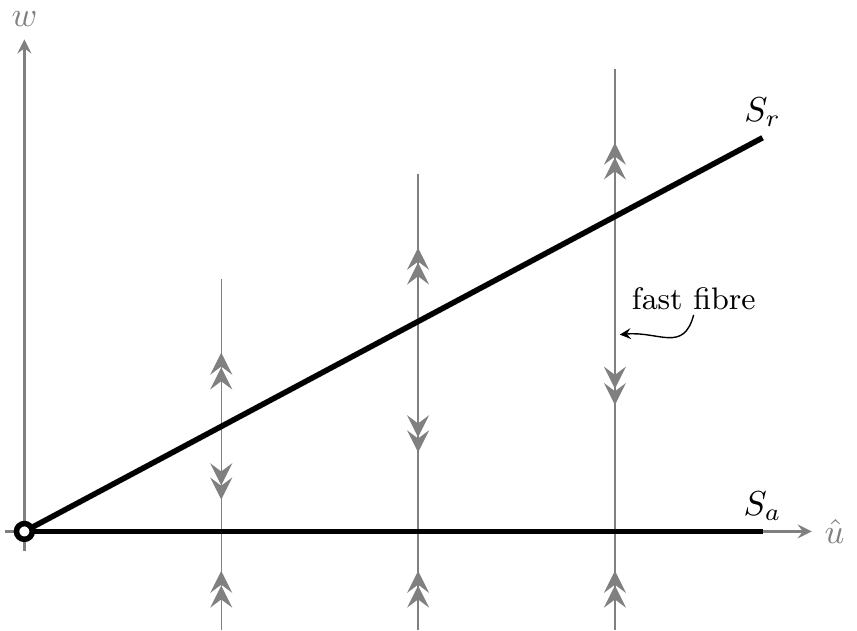}
\includegraphics[width=0.49\textwidth,valign=t]{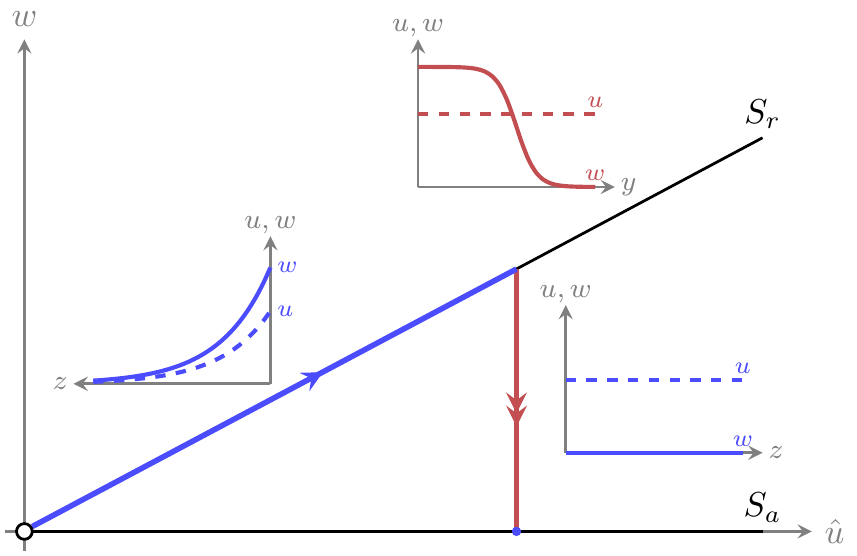}
\caption{The critical manifold $S$, projected into $(\tilde{u},w)$-space and the evolution of the original, fast variables $(u,w)$ in the different regions.  The open circle at the origin signifies that the original system \eref{eq:KS} has a removable singularity at this point, see Remark~\ref{rem:origin}.}
\label{fig:man-schematic}
\end{figure}

The first equation of \eref{eq:layer} gives $u = \tilde{u}/c + \alpha{\rm e}^{-cy/\mu}$.  However, as $y \to \pm\infty$ we require $u \to \tilde{u}/c$ and hence, $\alpha = 0$.  Therefore, along a fast fibre we have that $u = \tilde{u}/c$, while the evolution of $v$ and $w$ along the fast fibres is described by the second and third equation of \eref{eq:layer}.  An illustration is given in the left-hand panel of Figure~\ref{fig:man-schematic}.

\subsection{Reduced problem}
\label{subsec:reduced}

The three algebraic constraints of \eref{eq:reduced} are equivalent to the steady states of \eref{eq:layer}.  Consequently, the flow of the reduced problem is restricted to $S$.  We consider the flow on the two branches separately.  Firstly, on $S_a$ we have $\tilde{u}_z = 0$.  Therefore, there is no flow along $S_a$ and, using the asymptotic boundary conditions \eref{eq:bc}, we have $(u,v,w,\tilde{u}) = (u_r,0,0,cu_r)$.  Note that this also implies that $u = u_r$ and $\tilde{u} = cu_r$ along a fast fibre.  

Secondly, on $S_r$ we have $\tilde{u}_z = c\tilde{u}/\chi$, which can be solved exactly to give
\[ \tilde{u} = {\rm e}^{c(z + z^{\ast})/\chi}, \]
where $z^{\ast}$ is the constant of integration.  Consequently, 
\[ u = \frac{1}{c}{\rm e}^{c(z + z^{\ast})/\chi}, \quad v = \frac{1}{\chi}{\rm e}^{c(z + z^{\ast})/\chi}, \quad w = \frac{c}{\chi K}{\rm e}^{c(z + z^{\ast})/\chi}.  \]
We are free to choose $z^{\ast}$ since the problem is translation invariant.  To be consistent with the exact solution \eref{eq:limit}, we take $z^{\ast} = \chi\ln{(cu_r)}/c$.  Thus, in terms of the original variables $u$ and $w$, in the singular limit $\varepsilon \to 0$ the slow flow is described by 
\begin{equation}\label{eq:reduced-sol} 
u(z) = \begin{cases} u_r{\rm e}^{cz/\chi} & \textrm{on } S_r, \\ u_r & \textrm{on } S_a, \end{cases} \quad {\rm and} \quad
w(z) = \begin{cases} \dfrac{c^2 u_r}{K\chi}{\rm e}^{cz/\chi} & \textrm{on } S_r, \\ 0 & \textrm{on } S_a. \end{cases} 
\end{equation}
This coincides with \eref{eq:limit}, with the transition from $S_r$ to $S_a$ occurring at $z = 0$.  

\subsection{Singular heteroclinic orbits}

We now have enough information to construct heteroclinic orbits in the singular limit $\varepsilon \to 0$.  These {\em singular} orbits are concatenations of components from the reduced and layer problems.  Since the end state $u_r$ is a free parameter, we construct the waves in backward $z$.  

In backward $z$, a solution begins on $S_a$ from a point $(u,v,w,\tilde{u}) = (u_r,0,0,cu_r)$.  Since there is no evolution of the slow variables on $S_a$, the only possibility is for the solution to switch onto a fast fibre of the layer problem.  This connects the solution to the appropriate point on $S_r$: $(u,v,w,\tilde{u}) = (u_r,cu_r/\chi,c^2 u_r/(\chi K),cu_r)$.  Once back on $S_r$, the slow flow of the reduced problem evolves the solution towards the initial state of the wave $(u,v,w,\tilde{u}) = (0,0,0,0)$.  See the right-hand panel of Figure~\ref{fig:man-schematic} for an illustration.   

\subsection{Heteroclinic orbits for $0 < \varepsilon \ll 1$}
\label{subsec:persistence}

The persistence of the singular heteroclinic orbits for sufficiently small $0 < \varepsilon \ll 1$ is guaranteed by Fenichel theory \cite{Fenichel_71, Fenichel_79}.  Firstly, we consider the slow segments of the solutions.  Since $S_r$ and $S_a$ are normally hyperbolic, they deform smoothly to $\mathcal{O}(\varepsilon)$ close, locally invariant manifolds $S_{r,\varepsilon}$ and $S_{a,\varepsilon}$.  In this case, the model is simple enough that we can compute these manifolds explicitly, to any order:
\[ \begin{aligned}
S_{r,\varepsilon} &= \left\{ (u_{\varepsilon},v_{\varepsilon},w_{\varepsilon},\tilde{u}) \, \middle| \, u_{\varepsilon} = \frac{\tilde{u}(\chi - \varepsilon)}{c(\chi - \varepsilon(1 - \mu))}, v_{\varepsilon} = \frac{\tilde{u}}{\chi - \varepsilon(1 - \mu)}, w_{\varepsilon} = \frac{c\tilde{u}}{K (\chi - \varepsilon)} \right\}, \\
S_{a,\varepsilon} &= \left\{ (u_{\varepsilon},v_{\varepsilon},w_{\varepsilon},\tilde{u}) \, \middle| \, u_{\varepsilon} = \frac{\tilde{u}}{c}, v_{\varepsilon} = 0, w_{\varepsilon} = 0 \right\} = S_a.
\end{aligned} \]

It is not surprising that $S_{a,\varepsilon} = S_a$, since $S_a$ coincides with the background states of \eref{eq:KS}, which are not affected by the size of $\varepsilon$.  Consequently, the flow on $S_{a,\varepsilon}$ also remains unchanged, that is, there is no flow along $S_{a,\varepsilon}$.  On the other hand, the flow on $S_{r,\varepsilon}$ will be an $\mathcal{O}(\varepsilon)$ perturbation of the flow on $S_r$.  Since $S_{r,\varepsilon} \to (0,0,0,0)$ as $\tilde{u} \to 0$, the solution evolving on $S_{r,\varepsilon}$ will still connect (in backward $z$) to the initial state of the perturbed wave.  

We now consider the fast segment of the solutions.  Once again by Fenichel theory, we know that the unstable manifold of $S_r$, $\mathcal{W}^U(S_r)$, perturbs smoothly for $0 < \varepsilon \ll 1$ to the nearby local unstable manifold $\mathcal{W}^U(S_{r,\varepsilon})$.  Similar is true for the stable manifold of $S_a$.  Furthermore, since the intersection between $\mathcal{W}^U(S_r)$ and $\mathcal{W}^S(S_a)$ is transverse, it will persist for $0 < \varepsilon \ll 1$ and hence the fast fibres persist, connecting points on $S_{r,\varepsilon}$ to points on $S_{a,\varepsilon}$.  

Therefore, the solution constructed in the singular limit persists as a nearby solution of \eref{eq:KS} for $D_u = \mu\varepsilon$, $D_w = \varepsilon$, with $\varepsilon$ sufficiently small.  However, note that since $S_{a,\varepsilon}$ corresponds to a line of fixed points, the perturbed wave will connect to an end state $u_r(\varepsilon)$, $\mathcal{O}(\varepsilon)$ close to the original end state $u_r$ of the unperturbed wave.  Alternatively, since $u_r$ is likely to be a fixed quantity, we can say that the perturbed wave connects the original end states of the unperturbed wave but with a different speed $c(\varepsilon)$, $\mathcal{O}(\varepsilon)$ close to the original speed $c$.  
\qed\end{proof}

\begin{remark}
It is {\em a priori} not clear that {\sc gspt} extends to the singular point $(0,0)$.  However, using the methods of \cite{Doelman_Gardner_Kaper_01}, in which the the authors study a generalised Gierer--Meinhardt equation with a similar singularity, it can be shown that the theory indeed extends.  We refrain from going into the details.  
\end{remark}

\begin{remark}
The above results hold for $\mu = 0$.  Moreover, in this case we can solve the layer problem explicitly:
\[ u = \frac{\tilde{u}}{c}, \quad v = \frac{K\tilde{u}}{\chi K + \beta{\rm e}^{cy}}, \quad w = \frac{c\tilde{u}}{\chi K + \beta{\rm e}^{cy}}, \]
where $\beta$ is the integration constant.  
\end{remark}

\section{Conclusion}

Using {\sc gspt}, we proved the existence of travelling wave solutions to \eref{eq:KS} with $D_u = \mu\varepsilon$, $D_w = \varepsilon$ and $\varepsilon$ sufficiently small.  To leading order these solutions are given by \eref{eq:reduced-sol}, which are equivalent to the exact solutions of \cite{Feltham_Chaplain_00} given in \eref{eq:limit}.  This demonstrates the power of {\sc gspt} for studying the existence of travelling wave solutions to models such as the Keller--Segel model, even if exact solutions are not known.   

\subsection*{Acknowledgements}  This research was partially supported by the Australian Research Council's Discovery Projects funding scheme (project number DP110102775).  

\bibliographystyle{plain}
\bibliography{papers}

\end{document}